\documentclass[11pt]{amsart}
\usepackage{mathrsfs}
\usepackage{amssymb}

\usepackage{titletoc}
\input xy
  \xyoption{all}
\usepackage{amscd}
\usepackage{amsmath, amssymb}
\usepackage{amsfonts}
\usepackage[colorlinks,linkcolor=blue,citecolor=red, pdfstartview=FitH]{hyperref}

\numberwithin{equation}{section}

\theoremstyle{plain}
\newtheorem{thm}{Theorem}[section]
\newtheorem{theorem}[thm]{Theorem}
\newtheorem{lemma}[thm]{Lemma}
\newtheorem{corollary}[thm]{Corollary}
\newtheorem{proposition}[thm]{Proposition}

\theoremstyle{definition}

\newtheorem{question}[thm]{Question}

\newtheorem{remark}[thm]{Remark}

\newtheorem{definition}[thm]{Definition}

\newtheorem{example}[thm]{Example}

\newtheorem{defn-thm}[thm]{Definition-Theorem}

\newcommand{\sA}{{\mathcal A}}
\newcommand{\sB}{{\mathcal B}}

\newcommand{\ssR}{{\mathfrak R}}

\newcommand{\C}{{\mathbb C}}

\newcommand{\N}{{\mathbb N}}

\newcommand{\R}{{\mathbb R}}

\renewcommand{\S}{{\mathbb S}}

\newcommand{\Z}{{\mathbb Z}}

\newcommand{\qtq}[1]{\quad\mbox{#1}\quad}
\newcommand{\bp}{\bar{\partial}}
\newcommand{\Om}{\Omega}

\newcommand{\ts}{\otimes}

\newcommand{\st}{\stackrel}
\newcommand{\btheorem}{\begin{theorem}}
\newcommand{\etheorem}{\end{theorem}}
\newcommand{\bproposition}{\begin{proposition}}
\newcommand{\eproposition}{\end{proposition}}
\newcommand{\bdefinition}{\begin{definition}}
\newcommand{\edefinition}{\end{definition}}
\newcommand{\bcorollary}{\begin{corollary}}
\newcommand{\ecorollary}{\end{corollary}}
\newcommand{\bproof}{\begin{proof}}
\newcommand{\eproof}{\end{proof}}
\newcommand{\bremark}{\begin{remark}}
\newcommand{\eremark}{\end{remark}}
\newcommand{\eexample}{\end{example}}
\newcommand{\bexample}{\begin{example}}
\newcommand{\la}{\langle}
\newcommand{\elemma}{\end{lemma}}
\newcommand{\blemma}{\begin{lemma}}
\newcommand{\ra}{\rangle}
\newcommand{\sq}{\sqrt{-1}}

\newcommand{\om}{\omega}

\newcommand{\p}{\partial}

\renewcommand{\bar}{\overline}

\renewcommand{\phi}{\varphi}

\newcommand{\ee}{\end{eqnarray*}}
\newcommand{\be}{\begin{eqnarray*}}

\newcommand{\beq}{\begin{equation}}
\newcommand{\eeq}{\end{equation}}

\newcommand{\bd}{\begin{enumerate}}
\newcommand{\ed}{\end{enumerate}}

\renewcommand{\hat}{\widehat}
\renewcommand{\tilde}{\widetilde}

\renewcommand{\rm}{\textrm}

\newcommand{\pzj}{\frac{\p}{\p z^j}}
\newcommand{\pzi}{\frac{\p}{\p z^i}}
\newcommand{\pzk}{\frac{\p}{\p z^k}}
\newcommand{\pzl}{\frac{\p}{\p z^\ell}}
\newcommand{\pzs}{\frac{\p}{\p z^s}}
\newcommand{\pzt}{\frac{\p}{\p z^t}}

\newcommand{\bpzj}{\frac{\p}{\p \bar z^j}}
\newcommand{\bpzi}{\frac{\p}{\p \bar z^i}}
\newcommand{\bpzk}{\frac{\p}{\p \bar z^k}}
\newcommand{\bpzl}{\frac{\p}{\p \bar z^\ell}}

\newcommand{\bpzt}{\frac{\p}{\p \bar z^t}}

\renewcommand{\>}{\rightarrow}

\newcommand{\bi}{\bar{i}}
\newcommand{\bj}{\bar{j}}
\newcommand{\bk}{\bar{k}}

\newcommand{\Ga}{{\Gamma}}

\newcommand{\LC}{\nabla^\mathrm{LC}}

\renewcommand{\pzj}{\frac{\p}{\p z^j}}
\renewcommand{\pzi}{\frac{\p}{\p z^i}}
\renewcommand{\pzk}{\frac{\p}{\p z^k}}
\renewcommand{\pzl}{\frac{\p}{\p z^\ell}}
\renewcommand{\pzs}{\frac{\p}{\p z^s}}
\renewcommand{\pzt}{\frac{\p}{\p z^t}}

\renewcommand{\bpzj}{\frac{\p}{\p \bar z^j}}
\renewcommand{\bpzi}{\frac{\p}{\p \bar z^i}}
\renewcommand{\bpzk}{\frac{\p}{\p \bar z^k}}
\renewcommand{\bpzl}{\frac{\p}{\p \bar z^\ell}}

\renewcommand{\bpzt}{\frac{\p}{\p \bar z^t}}

\renewcommand{\>}{\rightarrow}

\renewcommand{\Ga}{\Gamma}

\newcommand{\Ta}{\Theta}
\newcommand{\ta}{\theta}
\newcommand{\La}{\Lambda}
\newcommand{\lam}{\lambda}

\newcommand{\w}{{\wedge}}
\newcommand{\bz}{{\bar{z}}}
\renewcommand{\bi}{{\bar{i}}}
\renewcommand{\bj}{{\bar{j}}}
\renewcommand{\bk}{{\bar{k}}}
\newcommand{\bel}{{\bar{\ell}}}
\renewcommand{\p}{{\partial}}
\renewcommand{\bp}{{\bar{\partial}}}
\newcommand{\bbp}{{\bar{p}}}
\newcommand{\bbq}{{\bar{q}}}

\newcommand{\bthm}{\begin{thm}}
\newcommand{\ethm}{\end{thm}}

\newcommand{\blem}{\begin{lem}}
\newcommand{\elem}{\end{lem}}

\newcommand{\brmk}{\begin{rmk}}
\newcommand{\ermk}{\end{rmk}}

\newcommand{\bpf}{\begin{proof}}
\newcommand{\epf}{\end{proof}}
\newcommand{\skpf}{\begin{proof}[Sketch of the proof]}

\newcommand{\I}{\sqrt{-1}}
\renewcommand{\lam}{\lambda}
\renewcommand{\La}{\Lambda}
\renewcommand{\st}{\star}

\renewcommand{\bz}{\bar{z}}
\newcommand{\bm}{\bar{m}}
\newcommand{\bn}{\bar{n}}

\newcommand{\cd}{\nabla}

\newcommand{\half}{\frac{1}{2}}

\newcommand{\dt}{\nabla^t}

\renewcommand{\cd}{\nabla}
\newcommand{\ssRic}{{\mathfrak{Ric}}}
\renewcommand{\st}{*}

\usepackage{fancyhdr}
\begin{document}
\title{ Curvatures of real connections on Hermitian manifolds} \makeatletter
\let\uppercasenonmath\@gobble
\let\MakeUppercase\relax
\let\scshape\relax
\makeatother

\author{Jun Wang}
\date{}
\address{Academy of Mathematics and Systems Science, Chinese Academy of Sciences, Beijing,
100190, China.}
\author{Xiaokui Yang}
\date{}
\address{Department of Mathematics and Yau Mathematical Sciences Center, Tsinghua University,  Beijing 100084,
China}
\email{\href{mailto:xkyang@mail.tsinghua.edu.cn}{{xkyang@mail.tsinghua.edu.cn}}}

\maketitle

\begin{abstract} Let $(M,g,J)$ be a Riemannian manifold with a compatible integrable complex structure $J\in\mathrm{End}(T_\R M)$ and  $\sA_{g,J}$ be the space of real connections on $T_\R M$ preserving both $g$ and
$J$. In this paper, we investigate the relationship between the
geometry of real connections in $\sA_{g,J}$ and that of Hermitian
connections on $T^{1,0}M$. In particular, we study the geometry of
real Chern connection $\nabla^{\mathrm{Ch,\R}}$ on $(M,g,J)$,  and
obtain K\"ahler-Einstein metrics by using real Chern-Einstein
metrics.
\end{abstract}

{\small\setcounter{tocdepth}{1} \tableofcontents}

\section{Introduction}

Let $(M,g,J)$ be a Riemannian manifold with a compatible integrable
complex structure $J\in \mathrm{End}(T_\R M)$ and $\sA_g$ be the
space of real connections on $T_\R M$ compatible with $g$. Let $h$
be the corresponding Hermitian metric of $(g,J)$ and $\sB_h$ be the
space of affine connections on the holomorphic tangent bundle
$T^{1,0}M$ compatible with $h$. For any $\nabla\in \sA_g$, we can
extend it to $T_\C M$ in the $\C$-linear way. The restriction of the
complexified connection $\nabla$ to the holomorphic tangent bundle
$T^{1,0}M$ is denoted by $\hat \nabla$. It is obvious that
$\hat\nabla\in \sB_h$. This gives a natural projection
$\pi:\sA_g\>\sB_h$ and it is easy to see that \beq \sA_g\cong
\sB_h\times \Gamma\left(M,\Om^{1,0}(\mathrm{Hom}\left(T^{1,0}M,
T^{0,1}M\right)\right).\eeq

\noindent Let $\sA_{g,J}$ be the space of real connections on $T_\R
M$ compatible with both $g$ and $J$. One can see that there is an
isomorphism $\rho:\sA_{g,J}\>\sB_h$ and \beq \sA_{g,J}\cong \sB_h
\cong\Gamma\left(M,\Om^{1,0}\left(\mathrm{End}\left(T^{1,0}M\right)\right)\right).\label{correspondence}\eeq

\noindent In the field of complex geometry,  several classical
connections in $\sB_h$ are investigated extensively in literatures.
For instances, \bd
\item The Chern connection $\nabla^{\mathrm{Ch}}$:  the
unique connection compatible with the Hermitian metric $h$ and also
the holomorphic structure $\bp$.

\item The Strominger-Bismut connection $\nabla^{\mathrm{SB}}$ (\cite{Str86, Bis89}).

\item $\hat\nabla^{\mathrm{LC}}$, the restriction of the
complexified Levi-Civita connection $\nabla^{\mathrm{LC}}$ to
$T^{1,0}M$. \ed

\noindent When $(M,g,J)$ is K\"ahler, all these connections are the
same.\\

  It is well-known that $\LC\in \sA_{g,J}$ if and only if
$\nabla^{\mathrm{LC}}J=0$, i.e. $(M,g,J)$ is  K\"ahler. Although
$\rho:\sA_{g,J}\>\sB_h$ is an isomorphism, the relationship between
the geometry of real connections in $\sA_{g,J}$ and that of
Hermitian connections in $\sB_h$ is still mysterious. For instance,
we set \beq
\nabla^{\mathrm{Ch,\R}}:=\rho^{-1}(\nabla^{\mathrm{Ch}})\in
\sA_{g,J}\eeq and it is called the \textbf{real Chern connection}.
It is obvious that $\nabla^{\mathrm{Ch,\R}}\neq \LC$ when $(M,g,J)$
is not K\"ahler. For any $X,Y,Z,W\in T_\R M$, the curvature of
$\nabla^{\mathrm{Ch,\R}}$ is  $$
R^{\mathrm{Ch,\R}}(X,Y,Z,W)=g\left(\nabla^{\mathrm{Ch,\R}}_X\nabla^{\mathrm{Ch,\R}}_YZ-\nabla^{\mathrm{Ch,\R}}_Y\nabla^{\mathrm{Ch,\R}}_XZ-\nabla^{\mathrm{Ch,\R}}_{[X,Y]}Z,W\right).
$$ The  \textbf{real Chern-Ricci curvature} $\mathrm{Ric}(\nabla^{\mathrm{Ch,\R}},g)$ is defined
by using the Riemannian metric $g$, i.e.
$$\mathrm{Ric}(\nabla^{\mathrm{Ch,\R}},g)(X,Y)=\sum_{i=1}^{2n} R^{\mathrm{Ch,\R}}(X,e_i,e_i,Y)$$
where $\{e_i\}_{i=1}^{2n}$ is an orthonormal frame with respect to
$g$.
 \bdefinition
$(M,g,J,\nabla^{\mathrm{Ch,\R}})$ is called  \textbf{real
Chern-Einstein} if \beq
\mathrm{Ric}(\nabla^{\mathrm{Ch,\R}},g)=\lambda g \qtq{for some
$\lambda\in \R.$}\eeq  If
$\mathrm{Ric}(\nabla^{\mathrm{Ch,\R}},g)=0$, it is also called
\textbf{real Chern-Ricci flat}. Moreover,
$(\nabla^{\mathrm{Ch,\R}},g)$ has positive  real Chern-Ricci
curvature if
$\mathrm{Ric}(\nabla^{\mathrm{Ch,\R}},g)\in\Gamma(M,\mathrm{Sym}^{\ts
2}T_\R M)$ and it is positive  definite. The negativity can be
defined similarly. \edefinition

\noindent 
There is a natural question concerning the relationship between real
Chern-Einstein metrics and K\"ahler-Einstein metrics.

\begin{question}\label{general} Let $(M,g,J)$ be a compact Hermitian manifold.
Suppose $(\nabla^{\mathrm{Ch,\R}}, g)$ is real Chern-Einstein, is
necessarily $(M,g,J)$ K\"ahler-Einstein? More generally, if
$\mathrm{Ric}(\nabla^{\mathrm{Ch,\R}}, g)$ is positive (resp.
negative), is $c_1(M)>0$ (resp. $c_1(M)<0$)?\end{question}

\noindent We can also ask similar questions for other connections in
$\sA_{g,J}$ or  $\sA_g$.  One shall see that the answer to the above
question is quite involved. Moreover, the existence of real
Chern-Ricci flat metrics is significantly different from others.

\btheorem\label{main1} Let $(M,g,J)$ be a compact Hermitian
manifold. Suppose $(\nabla^{\mathrm{Ch,\R}}, g)$ is real
Chern-Einstein with constant $\lambda\in \R$. If $\lambda\neq 0$,
then $(M,g,J)$ is K\"ahler-Einstein. \etheorem

\btheorem\label{main0} Let $(M,g,J)$ be a compact Hermitian
manifold. Suppose $\mathrm{Ric}(\nabla^{\mathrm{Ch,\R}},g)$ is
positive (resp. negative), then $M$ is projective and $c_1(M)>0$
(resp. $c_1(M)<0$).

\etheorem

\noindent However, when $(\nabla^{\mathrm{Ch,\R}}, g)$ is real
Chern-Ricci flat, i.e. $\lambda=0$,  $(M,g,J)$ is not necessarily
K\"ahler-Ricci flat.  We  construct explicit real Chern-Ricci flat
metrics on $\S^{2n-1}\times \S^1$. On the contrary, it is well-known
that, there is \emph{no} Chern-Ricci flat metric on $\S^{2n-1}\times
\S^1$ since $c_1^{\mathrm{BC}}(\S^{2n-1}\times \S^1)$ is not zero.

\btheorem\label{main4} There exist real Chern-Ricci flat metrics on
Hopf manifold $\S^{2n-1}\times \S^1$. \etheorem

\noindent We also obtain a classification for manifolds with real
Chern-Ricci flat metrics in Theorem \ref{main7}. Indeed, when
$\dim_\C M=2$, the Hopf surface is the only non-K\"ahler surface
which can support real Cher-Ricci flat metrics.\\

Question \ref{general} can also be proposed for arbitrary $\nabla\in
\sA_{g,J}$. Indeed, for any $\nabla\in \sA_{g,J}$, there exists some
$A\in\Gamma(M,\Om^1(\mathrm{End}(T_\R M)))$ such that
$$\nabla=\nabla^{\mathrm{Ch,\R}}+A,$$
and we also set $\hat\nabla=\rho(\nabla)\in \sB_h$. One can deduce
from (\ref{correspondence}) that $\nabla$ is determined by some
$\theta\in
\Gamma\left(M,\Om^{1,0}\left(\mathrm{End}\left(T^{1,0}M\right)\right)\right)$.
We write $\nabla^\theta$ for $\nabla\in\sA_{g,J}$, and
$\hat\nabla^\theta$ for $\hat\nabla\in\sB_h$. The
 curvature tensors of $\nabla^\theta$ and $\hat\nabla^\theta$ are denoted by $R^\theta,
 \mathfrak{R}^\theta$. In local holomorphic coordinates
$\{z^i\}$ of $M$, $\hat\nabla^\theta=\rho(\nabla^\theta)\in \sB_h$
is given by \be \hat \nabla^{\theta}_{\frac{\p}{\p z^i}}
\pzj=\left(\Gamma_{ij}^k+\theta_{ij}^k\right)\pzk,\ \ \ \ \
\hat\nabla^{\theta}_{\frac{\p}{\p \bar z^i}}\pzj=-h_{j\bar
q}h^{k\bar p}\bar{\theta_{ip}^q}\pzk \ee where
$\theta=\theta_{ij}^kdz^i\ts dz^j\ts{\pzk}$ and
$\Gamma_{ij}^k=h^{k\bar\ell}\frac{\p h_{j\bar\ell}}{\p z^i}$ is the
Christoffel symbol of the Chern connection. We also use  conventions
$R^\theta_{i\bar j k\bar\ell}$ and $\ssR^\theta_{i\bar j k\bar
\ell}$ for the components of $R^\theta$ and $\ssR^\theta$
respectively. The curvature of the Chern connection
$\nabla^{\mathrm{Ch}}$ is denoted by $\Theta$. We set
$\ssRic^{(1)}(\theta)=\sq \left(h^{k\bar\ell}\ssR^\theta_{i\bar
jk\bar\ell}\right)dz^i\wedge d\bar z^j$ and similarly, we  denote
the first Chern-Ricci curvature of the Chern connection
$\nabla^{\mathrm{Ch}}$ by $\Theta^{(1)}$.

\btheorem\label{main2} For any $\nabla^\theta\in \sA_{g,J}$,
curvature tensors $R^\theta$ and $\ssR^\theta$ are determined by
\be\begin{cases} R^\theta_{i\bar jk\bar\ell}&=\ssR^\theta_{i\bar j
k\bar\ell} =\Ta_{i\bj
k\bel}-\left(h_{k\bbp}\frac{\p\bar{\ta_{j\ell}^p}}{\p
z^i}+h_{p\bel}\frac{\p\ta_{ik}^p}{\p
\bz^{j}}\right)+\left(\ta_{ik}^p\bar{\ta_{j\ell}^q}h_{p\bbq}-h^{m\bn}\ta_{im}^{p}\bar{\ta_{jn}^q}h_{p\bel}h_{k\bbq}\right),\\
R_{ijk\bar\ell}^\theta&=\ssR^\theta_{ijk\bel}=\left(\frac{\p\ta_{jk}^m}{\p
z^i}-\frac{\p\ta_{ik}^m}{\p
z^j}+\Gamma_{jk}^s\theta_{is}^{\ell}-\Gamma_{js}^{\ell}\theta_{ik}^s+\Gamma_{is}^{\ell}\theta_{jk}^s
-\Gamma_{ik}^s\theta_{js}^{\ell}\right)h_{m\bel}\\&+\left(\ta_{ip}^{m}\ta_{jk}^p-\ta_{ik}^p\ta_{jp}^{m}\right)h_{m\bel}.
 \end{cases}\ee
Moreover, we have
 $[\ssRic^{(1)}(\theta)]=[\Theta^{(1)}]\in
H^{1,1}_{\mathrm{AC}}(M,\R)$ and \be
\ssRic^{(1)}(\ta)=\Ta^{(1)}-\I(\p\bar{\ta}_1-\bp\ta_1) \ee where
$\ta_1=\ta_{ik}^kdz^i$. \etheorem

  By using the isomorphism  $\rho:\sA_{g,J}\>\sB_h$,  there exists a \emph{unique} linear  family in $\{\nabla^t\}_{t\in \R}\subset \sA_{g,J}$
 which connects the real Chern connection $\nabla^{\mathrm{Ch,\R}}=\rho^{-1}(\nabla^{\mathrm{Ch}})$
 and the real Strominger-Bismut connection
 $\nabla^{\mathrm{SB,\R}}=\rho^{-1}(\nabla^{\mathrm{SB}})$, and it is given by
\beq
\nabla^t=(1-t)\nabla^{\mathrm{Ch,\R}}+t\nabla^{\mathrm{SB,\R}}.\eeq
This family is  introduced in \cite{Gau97} by  Gauduchon, and we
call it \emph{Gauduchon connections}. They are systematically
investigated in recent papers \cite{FY08, FTY09, AGF12, BM13, AGF14,
FY15,OUV17,AOUV18, FZ19, YZZ19, ZZ19a} on the construction of
invariant solutions to the Strominger systems on complex Lie groups
and their quotients. A straightforward computation shows \beq
\rho\left(\nabla^0\right)=\nabla^{\mathrm{Ch}},\ \
\rho\left(\nabla^{\frac{1}{2}}\right)=\hat\nabla^{\mathrm{LC}},\ \
\rho\left(\nabla^1\right)=\nabla^{\mathrm{SB}}.\eeq Hence, these
classical connections are all in the Gauduchon family. The curvature
relations for $\nabla^{\mathrm{LC}}, \nabla^{\mathrm{SB}}$ and
$\hat\nabla^{\mathrm{LC}}$ are extensively investigated and have
been formulated in differential notions (e.g. \cite{Yau74,
Gra76,Gau77a, Gau77b, TV81, Gau84, AD99,Fu12, LY12,LY16,WYZ16,
Yang17a,AOUV18, HLY18,YZ18a, YZ18b}). We shall formulate them in a
uniform way for readers' convenience. As usual, the curvature tensor
of $\rho(\nabla^t)\in\sB_h$ is denoted by $\ssR(\omega,t)$ and some
notions can be found in \cite{LY16}.

\bcorollary\label{god}
 The curvature tensor of Gauduchon connection $\rho(\nabla^t)\in\sB_h$ is
$$ \ssR_{i\bj k\bel}(\omega, t) =\Ta_{i\bj k\bel}+t(\Ta_{i\bel k\bj}+\Ta_{k\bj
i\bel}-2\Ta_{i\bj k\bel})+t^2(T_{ik}^p\bar{T_{j\ell
}^q}h_{p\bbq}-h^{p\bbq}h_{m\bar\ell} h_{k\bar
n}T_{ip}^m\bar{T_{jq}^n})$$ where
$T_{ij}^k=\Gamma_{ij}^k-\Gamma_{ji}^k$ are the components of the
torsion tensor of $\nabla^{\mathrm{Ch}}$. The Ricci curvatures are
given by
 \beq
 \ssRic^{(1)}(\omega,t)=\Ta^{(1)}-t(\p\p^\st\om+\bp\bp^\st\om).\label{FZ}
 \eeq
\noindent and \be
\begin{cases}
\ssRic^{(2)}(\omega,t)=&\Ta^{(1)}-(1-2t)\I
\La\p\bp\om-(1-t)(\p\p^\st\om+\bp\bp^\st\om)
\\&+(1-t)^2\I T\boxdot T-t^2\I T\circ T,\\
 \ssRic^{(3)}(\omega, t)=&\Ta^{(1)}-t\I
\La\p\bp\om-(1-t)\p\p^\st\om-t\bp\bp^\st\om
\\&+(t-t^2)\I T\boxdot T+t^2 T((\p^\st\om)^{\#}),\\
\ssRic^{(4)}(\omega, t)=&\Ta^{(1)}-t\I
\La\p\bp\om-(1-t)\bp\bp^\st\om-t\p\p^\st\om
\\&+(t-t^2)\I T\boxdot T+t^2\bar{ T((\p^\st\om)^{\#})}.
\end{cases}\label{examplecorrespondence}
\ee

\noindent The scalar curvatures are related by
 \beq \label{FZ2}\begin{cases}
  s^{(1)}(\omega, t)=s_{\mathrm{C}}-2t\la \p\p^\st\om,\om\ra;\\
 s^{(2)}(\omega, t)=s_{\mathrm{C}}-(1-2t)\la
 \p\p^\st\om,\om\ra-t^2(2|\p\om|^2+|\p^\st\om|^2),
\end{cases}
 \eeq
where we use the following notations for scalar curvatures
$$s^{(1)}(\omega,t)=h^{i\bar j}h^{k\bar\ell}\ssR_{i\bar
jk\bar\ell}(\omega,t),\ \ s^{(2)}(\omega,t)=h^{i\bar\ell}h^{k\bar
j}\ssR_{i\bar j k\bar\ell}(\omega,t).$$ \ecorollary

\noindent
 Let $s$ be the Riemannian scalar curvature of the background Riemannian manifold $(M,g)$.  One can deduce that (e.g. \cite[Corollary~3.7]{Yang20})
\beq s=2s_{\mathrm{C}}-2\sq\p^*\bp^*\omega-\frac{1}{2}|T|^2.\eeq By
using this formula and (\ref{examplecorrespondence}), one can get
scalar curvature relations for connections
$\nabla^{\mathrm{Ch}},\hat\nabla^{\mathrm{Ch}}$ and
$\nabla^\mathrm{SB}$ simultaneously.

\bremark We should point out that, many parts of Corollary \ref{god}
are known in literatures. Indeed,  in a recent paper \cite{FZ19},
J.-X. Fu and X.-C. Zhou established a variety of scalar curvature
formulas for Gauduchon connections for general almost Hermitian
manifolds, and (\ref{FZ}), (\ref{FZ2}) are also obtained in
\cite{FZ19}. For more discussions on scalar curvatures of almost
Hermitian manifolds, we refer to
\cite{AD99,Li10,Zha12,LU17,CZ18,FZ19} and the references therein.
One can also formulate curvatures relations in the conformal
setting(e.g. \cite{CRS14, ACS15, LM18, CCN19}).

 \eremark

As applications of curvature relations discussed above and methods
developed in \cite{LY16,Yang16, LY17,HLY18,Yang19,Yang19b, Yang20},
we obtain \btheorem\label{main3}
 Let $(M,\omega)$ be a compact Hermitian manifold and $\nabla^t$ be the Gauduchon connection on $M$. If $s^{(1)}(\omega, t)\geq 0$ for some $t>0$, then either
 \begin{enumerate}
 \item $\kappa(M)=-\infty $; or
 \item $ \kappa(M)=0$ and $(M,\om)$ is conformally balanced with $K_M$ is holomorphic torsion, i.e. $K_M^{\otimes m}=0$ for some $m\in \N_+$.
 \end{enumerate}
 \etheorem

\noindent By using Theorem \ref{main3}, we can classify
$t$-Gauduchon-Ricci flat surfaces. Recall that, a Hermitian manifold
$(M,\omega)$ is called \textbf{$t$-Gauduchon-Ricci flat} if \beq
\ssRic^{(1)}(\omega,t)=0 \qtq{for some} t\in \R.\eeq  When
$t=0,\frac{1}{2}$ and $1$, it is also called \emph{Chern-Ricci
flat}, \emph{Levi-Civita-Ricci flat} and
\emph{Strominger-Bismut-Ricci flat} respectively.

\btheorem\label{main5}  Let $S$ be a compact complex surface. If it
admits a $t$-Gauchon-Ricci flat metric $\omega$ for some $t>0$, then
$S$ is a minimal surface lying in one of the following:
\begin{enumerate}
\item Enriques surfaces.
\item bi-elliptic surfaces.
\item K3-surfaces.
\item  tori.
\item Hopf surfaces.
\end{enumerate}
\etheorem

\noindent Note that $(1)$-$(4)$ are all K\"ahler Calabi-Yau
surfaces. We also construct an explicit family of
$t$-Gauduchon-Ricci flat metrics on diagonal Hopf surfaces.

\btheorem\label{main6}  Let $ \omega_0$ be the canonical metric on
standard
 Hopf manifold $\S^{2n-1}\times \S^1$ ($n\geq 2$). Then the Hermitian metric \beq
\omega_t=\omega_0+4\left(\frac{2(n-1)t}{n}-1\right)\cdot\sq
\p\bp\log |z|^2,\ \ t>0 \eeq is $t$-Gauduchon-Ricci flat. That is,
 $\ssRic^{(1)}(\omega_t, t)=0$ for each $t>0$.  \etheorem

\noindent When $t=\frac{1}{2}$, the Hermitian metric \beq
\omega_{\mathrm{LC}}=\omega_0-\frac{4}{n}\cdot \sq\p\bp\log|z|^2\eeq
  is exactly the \emph{Levi-Civita-Ricci flat} metric on $\S^{2n-1}\times \S^1$ constructed in
\cite[Theorem~6.2]{LY16} (see also \cite[Theorem~7.3]{LY17}).
 When $t=1$, the Hermitian metric
\beq \omega_{\mathrm{SB}}=\omega_0-\frac{4(n-2)}{n}\cdot \sq
\p\bp\log|z|^2\eeq
  is a \emph{Strominger-Bismut-Ricci flat} metric on $\S^{2n-1}\times
  \S^1$.\\ 

 The paper is organized as follows. In Section \ref{background}, we
 recall basic materials for readers' convenience and fix our
 notations. In Section \ref{relation}, the spaces of connections are
 discussed. The curvatures of real connections in $\sA_{g,J}$ are
 computed in Section \ref{curvatureofagj} and Theorem \ref{main2} is
 obtained. In Section \ref{cherneinstein}, we establish Theorem \ref{main1}, Theorem \ref{main0} and Theorem \ref{main4}.
 In Section \ref{classification} we classify compact complex surfaces with $t$-Gauduchon-Ricci flat metrics and prove Theorem \ref{main3} and Theorem \ref{main5}. We construct $t$-Gauduchon-Ricci flat metrics in Section
 \ref{hopfmanifold}, and  in the Appendix, we include a detailed computation for Corollary \ref{god}.\\

\textbf{Acknowledgements.} We are very grateful to Professor
 K.-F. Liu and Professor S.-T. Yau for their support, encouragement and stimulating
discussions over  years. The second author would also like to thank
Teng Fei and Valentino Tosatti for helpful suggestions.

\vskip 2\baselineskip

\section{Background materials}\label{background}

In this section, we give some background materials for readers'
convenience.
\subsection{Levi-Civita connection and its complexification.} Let's recall some elementary settings.
Let $(M, g, \LC)$ be a $2n$-dimensional Riemannian manifold with the
Levi-Civita connection $\nabla^{\mathrm{LC}}$. The tangent bundle of
$M$ is also denoted by $T_\R M$. The Riemannian curvature tensor of
$(M,g,\LC)$ is
$$ R(X,Y,Z,W)=g\left(\LC_X\LC_YZ-\LC_Y\LC_XZ-\LC_{[X,Y]}Z,W\right)$$
for  $X,Y,Z,W\in T_\R M$. Let $T_\C M=T_\R M\ts \C$ be the
complexification. One can extend the metric $g$ and the Levi-Civita
connection $\LC$ to $T_{\C}M$ in the $\C$-linear way. Hence for any
$a,b,c,d\in \C$ and $X,Y,Z,W\in T_\C M$,  $$ R(aX,bY,cZ,
dW)=abcd\cdot R(X,Y,Z,W).$$

\noindent Let $(M,g,J)$ be an almost Hermitian manifold, i.e.,
$J:T_\R M\>T_\R M$ with $J^2=-1$, and
 for any $X,Y\in T_\R M$, $ g(JX,JY)=g(X,Y)$. The Nijenhuis tensor $N_J:\Gamma(M,T_\R M)\times \Gamma(M,T_\R
M)\>\Gamma(M,T_\R M)$ is defined as $$
N_J(X,Y)=[X,Y]+J[JX,Y]+J[X,JY]-[JX,JY].$$ The almost complex
structure $J$ is  \emph{integrable} if $N_J\equiv 0$ and then
$(M,g,J)$ is a Hermitian manifold. We  also extend $J$ to $T_\C M$
in the $\C$-linear way, i.e.  for any $X,Y\in T_\C M$, we still have
$ g(JX,JY)=g(X,Y).$
 By Newlander-Nirenberg's
theorem, there exist real coordinates $\{x^i,x^I\}$ such that
$z^i=x^i+\sq x^I$ are  local holomorphic coordinates on $M$.
 The Hermitian form $h:T_\C M\times T_\C M\>\C$ is given
by \beq h(X,Y):= g(X,Y),
 \ \ \ \ \ X,Y\in T_\C M.\label{complex}\eeq
By the $J$-invariant property of $g$, \beq
h_{ij}:=h\left(\frac{\p}{\p z^i},\frac{\p}{\p z^j}\right)=0,
\qtq{and} h_{\bar i\bar j}:=h\left(\frac{\p}{\p \bar
z^i},\frac{\p}{\p \bar z^j}\right)=0\eeq and \beq h_{i\bar
j}:=h\left(\frac{\p}{\p z^i},\frac{\p}{\p \bar
z^j}\right)=\frac{1}{2}\left(g_{ij}+\sq
g_{iJ}\right).\label{1003}\eeq It is obvious that $(h_{i\bar j})$ is
a positive Hermitian  matrix. Let $\omega$ be the fundamental
$2$-form associated to the $J$-invariant metric $g$,
 $ \omega(X,Y)=g(JX,Y)$.
   In local complex coordinates,
 $ \omega=\sq h_{i\bar j} dz^i\wedge d\bar z^j.$
We shall use the components of the complexified curvature tensor
$R$, for instances,  \beq R_{i\bar j k\bar
\ell}:=R\left(\frac{\p}{\p z^i}, \frac{\p}{\p \bar z^j},
\frac{\p}{\p z^k}, \frac{\p}{\p \bar z^\ell}\right),\ \ \ R_{i j k
\ell}:=R\left(\frac{\p}{\p z^i}, \frac{\p}{\p z^j}, \frac{\p}{\p
z^k}, \frac{\p}{\p z^\ell}\right).\ \ \ \eeq The components of
complexified curvature tensor have the same properties as the
components of the real curvature tensor, for instances: $ R_{i\bar j
k\ell}=-R_{\bar j i k\ell}$, $R_{i\bar j k\bar\ell}=R_{k\bar \ell
i\bar j}$ and in particular, the first Bianchi identity holds: $
R_{i\bar j k\bar \ell}+R_{ik\bar\ell \bar j}+R_{i\bar \ell \bar
jk}=0$.

\vskip 1\baselineskip

\subsection{The induced Levi-Civita connection on $(T^{1,0}M,h)$}
\noindent Since $T^{1,0}M$ is a subbundle of $T_{\C}M$, there is an
induced connection $\hat\nabla^{\mathrm{LC}}$ on $T^{1,0}M$ given by
\beq \hat\nabla^{\mathrm{LC}}=\pi\circ\nabla:
\Gamma(M,T^{1,0}M)\stackrel{\LC}{\rightarrow}\Gamma(M, T^*_{\C}M\ts
T_{\C}M)\stackrel{\pi}{\rightarrow}\Gamma(M,T^*_{\C}M\ts T^{1,0}M).
\label{metricconnection}\eeq Moreover, $\hat \nabla^{\mathrm{LC}}$
is a metric connection on the Hermitian holomorphic vector bundle
$(T^{1,0}M, h)$ and it is determined by the  relations \beq
\hat\nabla^{\mathrm{LC}}_{\frac{\p}{\p z^i}}\frac{\p}{\p
z^k}:=\hat\Gamma_{ik}^p\frac{\p}{\p z^p} \qtq{and}
\hat\nabla^{\mathrm{LC}}_{\frac{\p}{\p \bar z^j}}\frac{\p}{\p
z^k}:=\hat \Gamma_{\bar jk}^p\frac{\p}{\p z^p} \eeq where \beq \hat
\Gamma_{ij}^k=\frac{1}{2}h^{k\bar \ell}\left(\frac{\p h_{j\bar
\ell}}{\p z^i}+\frac{\p h_{i\bar \ell}}{\p z^j}\right), \qtq{and}
\hat \Gamma_{\bar i j}^k=\frac{1}{2} h^{k\bar \ell}\left( \frac{\p
h_{j\bar\ell}}{\p\bar z^i}-\frac{\p h_{j\bar i}}{\p\bar
z^\ell}\right).\eeq

\noindent The curvature  tensor $\mathfrak{R}\in \Gamma(M,\Lambda^2
T^*_{\C}M\ts T^{*1,0}M\ts T^{1,0}M)$ of $\hat\nabla^{\mathrm{LC}}$
is given by \beq \mathfrak{R}(X,Y)s
=\hat\nabla^{\mathrm{LC}}_{X}\hat\nabla^{\mathrm{LC}}_Ys-\hat\nabla^{\mathrm{LC}}_Y\hat\nabla^{\mathrm{LC}}_Xs-\hat\nabla^{\mathrm{LC}}_{[X,Y]}s\eeq
for any $X,Y\in T_{\C}M$ and $s\in T^{1,0}M$. The curvature tensor
$\mathfrak{R}$ has components
 \beq \begin{cases}\mathfrak{R}_{i\bar
jk}^{\ell}=-\left(\frac{\p \hat\Gamma^{\ell}_{ik}}{\p \bar
z^j}-\frac{\p \hat\Gamma^{\ell}_{\bar jk}}{\p z^i}+\hat\Gamma_{
ik}^{s}\hat\Gamma^{\ell}_{\bar js}-\hat\Gamma_{ \bar
jk}^{s}\hat\Gamma^{\ell}_{is}\right),\\\mathfrak{R}_{i
jk}^{\ell}=-\left(\frac{\p \hat\Gamma^{\ell}_{ ik}}{\p z^j}-\frac{\p
\hat\Gamma^{\ell}_{jk}}{\p
z^i}+\hat\Gamma_{ik}^{s}\hat\Gamma^{\ell}_{js}-\hat\Gamma_{
jk}^{s}\hat\Gamma^{\ell}_{is}\right),\\   \mathfrak{R}_{\bar i\bar j
k}^{\ell}=-\left(\frac{\p \hat\Gamma^{\ell}_{\bar ik}}{\p \bar
z^j}-\frac{\p \hat\Gamma^{\ell}_{\bar jk}}{\p \bar
z^i}+\hat\Gamma_{\bar ik}^{s}\hat\Gamma^{\ell}_{\bar
js}-\hat\Gamma_{\bar jk}^{s}\hat\Gamma^{\ell}_{\bar
is}\right).\end{cases}\eeq

\noindent With respect to the Hermitian metric $h$ on $T^{1,0}M$, we
use the  convention  \beq \mathfrak{R}_{\bullet\bullet k\bar
\ell}:=\sum_{s=1}^n\mathfrak{R}_{\bullet\bullet k}^{s}h_{s\bar
\ell}. \eeq

\bcorollary{\label{2nd}} One has the following relations
$$R_{ijk}^\ell=\mathfrak{R}_{ijk}^\ell,\ \ \ R_{\bar i\bar j
k}^{\ell}=\mathfrak{R}_{\bar i\bar j k}^\ell,$$ and \beq R_{i\bar j
k}^\ell=-\left(\frac{\p \hat \Gamma^{\ell}_{ik}}{\p \bar
z^j}-\frac{\p \hat\Gamma^{\ell}_{\bar jk}}{\p z^i}+\hat\Gamma_{
ik}^{s}\hat\Gamma^{\ell}_{\bar js}-\hat\Gamma_{ \bar
jk}^{s}\hat\Gamma^{\ell}_{s i}-\hat\Gamma_{\bar s
i}^\ell\cdot\bar{\hat\Gamma_{\bar k j}^s}\right)=\mathfrak{R}_{i\bar
j k}^\ell+\hat\Gamma_{i\bar s }^\ell\cdot\bar{\hat\Gamma_{j\bar k
}^s}.\eeq

\ecorollary

\vskip 1\baselineskip

\subsection{Curvature of the Chern connection on $(T^{1,0}M,h)$} On
the Hermitian holomorphic tangent bundle $(T^{1,0}M,h)$, the Chern
connection $\nabla^{\mathrm{Ch}}$ is the unique connection which is
compatible with the holomorphic structure and also the Hermitian
metric. The curvature tensor of $\nabla^{\mathrm{Ch}}$ is denoted by
$\Theta$ and it has components \beq \Theta_{i\bar j k\bar
\ell}=-\frac{\p^2 h_{k\bar \ell}}{\p z^i\p \bar z^j}+h^{p\bar
q}\frac{\p h_{p\bar \ell}}{\p \bar z^j}\frac{\p h_{k\bar q}}{\p
z^i}.\eeq
 It is
well-known that the \emph{(first) Chern Ricci curvature} \beq
\Theta^{(1)}:= \sq\Theta^{(1)}_{i\bar j} dz^i\wedge d\bar z^j\eeq
represents the first Chern class $ c_1(M)$ of $M$ where
  \beq
\Theta^{(1)}_{i\bar j}= h^{k\bar \ell}\Theta_{i\bar j k\bar \ell}
=-\frac{\p^2 \log \det(h_{k\bar \ell})}{\p z^i\p\bar
z^j}.\label{1chern}\eeq The \emph{second Chern Ricci curvature}
$\Theta^{(2)}=\sq \Theta^{(2)}_{i\bar j}dz^i\wedge d\bar z^j$ has
components \beq \Theta^{(2)}_{i\bar j}=h^{k\bar \ell}\Theta_{k\bar
\ell i\bar j}. \label{2chern}\eeq The \emph{Chern scalar curvature}
$s_\mathrm{C}$ of the Chern curvature tensor $\Theta$ is defined by
\beq s_{\mathrm{C}}=h^{i\bar j}h^{k\bar \ell}\Theta_{i\bar j k\bar
\ell}.\eeq Similarly, we can define \beq \Theta^{(3)}_{i\bar
j}=h^{k\bar\ell}\Theta_{i\bar \ell k \bar j},\ \ \
\Theta^{(4)}_{i\bar j}=h^{k\bar\ell}\Theta_{k\bar j
i\bar\ell}\qtq{and} s_{\mathrm{C}}^{(2)}=h^{i\bar
j}\Theta^{(3)}_{i\bar j}=h^{i\bar j}\Theta^{(4)}_{i\bar
j}.\label{34}\eeq

\noindent The following lemma is well-known (e.g.
\cite[Lemma~3.6]{Yang20}).
 \blemma Let $(X,\omega)$ be a compact
Hermitian manifold. Then \beq \la
\bp\bp^*\omega,\omega\ra=|\bp^*\omega|^2- \sq
\p^*\bp^*\omega.\label{c2} \eeq In particular, if $\omega$ is a
Gauduchon metric, we have \beq \la
\bp\bp^*\omega,\omega\ra=|\bp^*\omega|^2.\label{c3}\eeq \elemma

\noindent According to \cite[Theorem~4.1]{LY16}, one has

\btheorem\label{chernricci}  The  Chern-Ricci curvatures are related
by \beq
\begin{cases}
\Theta^{(2)}=\Theta^{(1)}-\sq\Lambda\left(\p\bp\omega\right)-(\p\p^*\omega+\bp\bp^*\omega)+{\sq}T\boxdot\bar
T,\\ \Theta^{(3)}=\Theta^{(1)}-\p\p^*\omega,\\
\Theta^{(4)}=\Theta^{(1)}-\bp\bp^*\omega\end{cases}\eeq where
$T\boxdot \bar T:=h^{p\bar q} h_{k \bar \ell}T_{ip}^k\cdot \bar
{T_{jq}^\ell} dz^i\wedge d\bar z^j$. The scalar curvatures are
related by
 \beq \begin{cases}
s^{(2)}_{\mathrm{C}}=s_{\mathrm{C}}-\la\p\p^*\omega, \omega\ra\\
s=2s_{\mathrm{C}}-2\sq\p^*\bp^*\omega-\frac{1}{2}|T|^2,\end{cases}\eeq
where $s$ is the scalar curvature of the background Riemannian
metric. \etheorem

%
%

\vskip 2\baselineskip
\section{The space $\sA_{g,J}$ of connections preserving  metrics and complex
structures}\label{relation}

\subsection{Spaces of real connections and Hermitian connections}
Let $(M,g, J)$ be a Riemannian manifold with  a compatible
integrable complex structure $J$.  We consider two spaces of affine
connections on the real tangent bundle $T_\R M$: \bd
\item[]
$\sA_{g}=\{\cd | \text{ $\cd$ is an affine connection on $T_\R M$
satisfying $\cd g=0$ } \}$,
\item[]
$\sA_{g,J}=\{\cd | \text{ $\cd$ is an affine connection on $T_\R M$
satisfying $\cd g=0$ and $\cd J=0$ } \}.$ \ed Clearly,
$\iota:\sA_{g,J}\hookrightarrow\sA_g$. In the following, we shall
give a characterization of connections in $\sA_{g,J}$.

\blemma Let $\cd^0\in \sA_{g,J}$. Suppose $\cd=\cd^{0}+A$ with
$A\in\Gamma(M,\Om^1(\mathrm{End}(T_\R M)))$. Then $\cd\in \sA_{g,J}$
if and only if \be\label{compatible-eqn-2} \left\{
\begin{aligned}
&[A,J]=0,
\\&g(AX,Y)+g(X,AY)=0 \quad \text{for any $X,Y\in \Gamma(M,T_\R{M})$}.
\end{aligned}
\right. \ee \elemma

\bpf It is easy to see that $$\cd J=\cd^{0} J+[A,J]=[A,J]$$ and  \be
g(\cd X,Y)+g(X,\cd Y)-d(g(X,Y))=g(AX,Y)+g(X,AY) \ee since
$\nabla^0J=0$ and $\nabla^0g=0$. \epf

Let $A_\C\in \Gamma(M,\Om^1(\mathrm{End}(T_\C M)))$ be the
complexification of $A\in\Gamma(M,\Om^1(\mathrm{End}(T_\R M)))$.
According to the decomposition $T_\C M =T^{1,0}M\bigoplus T^{0,1}M$,
$A_\C$ has a matrix representation $\left[\begin{array}{ll}
A_{11}&A_{12}\\
A_{21}&A_{22}
\end{array}\right]$, and $[A,J]=0$ implies $A_{12}=A_{21}=0$. Hence,
we have $A_\C=A_{11}+A_{22}$ where
$A_{11}\in\Gamma(M,\Om^1(\mathrm{End}(T^{1,0}M)))$ and $A_{22}\in
\Gamma(M,\Om^1(\mathrm{End}(T^{0,1}M)))$. Let
$A_{11}=\theta_1+\theta_2$ and $A_{22}=\theta_3+\theta_4$ where
$\theta_1\in\Gamma(M,\Om^{1,0}(\mathrm{End}(T^{1,0}M)))$,
$\theta_2\in\Gamma(M,\Om^{0,1}(\mathrm{End}(T^{1,0}M)))$,
 $\theta_3\in\Gamma(M,\Om^{1,0}(\mathrm{End}(T^{0,1}M)))$ and $\theta_4\in\Gamma(M,\Om^{0,1}(\mathrm{End}(T^{0,1}M)))$. Since $A$ is real and $\nabla h=0$, for any $X,Y\in
\Gamma(M,T^{1,0}M)$, we have \beq \begin{cases} h(\theta_{1}X,\bar
Y)+h(X,\theta_{3}\bar Y)=0,\\ h(\theta_{2}X,\bar
Y)+h(X,\theta_{4}\bar Y)=0,\\
\bar \theta_1=\theta_4,\\
\bar\theta_2=\theta_3.
\end{cases}\eeq
Hence, $A$ is determined by $\theta_1$. If we write
$\theta_1=\theta_{ij}^kdz^i\ts dz^j \ts \frac{\p}{\p z^k}$, then the
complexification of  $\nabla\in\sA_{g,J}$ is given by
\beq\label{realcompatible} \begin{cases}
{\nabla}_{\pzi}\pzj={\nabla^0}_{\pzi}\pzj+\ta_{ij}^k\pzk,
\\ {\nabla}_{\pzi}\bpzj={\nabla^0}_{\pzi}\bpzj-h_{q\bj}h^{p\bk}\ta_{ip}
^q\bpzk
\end{cases} \eeq and their conjugations. Therefore, one has

\bproposition $\sA_{g,J}\cong
\Gamma(M,\Om^{1,0}(\mathrm{End}(T^{1,0}M)))$. \eproposition

\noindent On the holomorphic tangent bundle $T^{1,0}M$, we consider
the following space: \beq \sB_{h}=\{\nabla |\nabla\ \text{is an
affine connection on $T^{1,0}M$ satisfying}\ \nabla h=0\}.\eeq It is
easy to see that the Chern connection $\nabla^{\mathrm{Ch}}\in
\sB_h$. Moreover, for any $\nabla\in \sA_{g,J}$, the restriction
$\hat\nabla:\Gamma(M,T^{1,0}M)\>\Gamma(M,\Om^1(T^{1,0}M))$ of its
complexification $\nabla:\Gamma(M,T_\C M)\>\Gamma(M,\Om^1(T_\C M))$
is in $\sB_h$. Hence, there is a natural map \beq \rho:
\sA_{g,J}\>\sB_h. \eeq

 \bcorollary $\rho: \sA_{g,J}\> \sB_h$ is an isomorphism. \ecorollary

\noindent Indeed, for any $\nabla\in \sB_h$, there exists some $B\in
\Gamma(M,\Om^{1}(\mathrm{End}(T^{1,0}M)))$ such that
$\nabla=\nabla^{\mathrm{Ch}}+B$. We have the decomposition
$B=B_{1}+B_2$ where $B_1=\theta_{ij}^k dz^i\ts dz^j\ts \frac{\p}{\p
z^k}\in\Gamma(M,\Om^{1,0}(\mathrm{End}(T^{1,0}M)))$ and
$B_2=\eta_{\bar i j}^k d\bar z^i\ts dz^j\ts \frac{\p}{\p z^k}\in
\Gamma(M,\Om^{0,1}(\mathrm{End}(T^{1,0}M)))$. Since $\nabla h=0$, we
deduce \beq\theta_{ij}^k h_{k\bar \ell}+\bar{\eta_{\bar i \ell}^p}
h_{j\bar p}=0.\eeq Hence, the connection $\nabla\in \sB_h$ is given
by \beq \begin{cases}
{\nabla}_{\pzi}\pzj={\nabla^{\mathrm{Ch}}}_{\pzi}\pzj+\ta_{ij}^k\pzk,
\\ {\nabla}_{\bpzi}\pzj=-h_{j\bar q}h^{k\bar p}\bar{\ta_{ip}
^q}\pzk.
\end{cases}
\eeq That means $\nabla\in \sB_h$ is determined by $B_1$, i.e.
$\sB_h\cong \Gamma(M,\Om^{1,0}(\mathrm{End}(T^{1,0}M)))$. By using
similar interpretations, one can show \bcorollary $\sA_g\cong
\sB_h\times \Gamma(M,\Om^{1,0}(\mathrm{Hom}(T^{1,0}M,T^{0,1}M)))$.
\ecorollary

\vskip 1\baselineskip

\subsection{The real correspondences for metric connections on
$T^{1,0}M$} On a Hermitian manifold $(M,g,J)$, we have the following
diagram for spaces of connections: \be \xymatrix{
\sA_{g,J}\ar[r]^{\iota}\ar[rd]_{\rho}&\sA_{g}\ar[d]^{\pi}
\\& \sB_h.
} \ee As we discussed in the previous section,
$\hat\nabla^{\mathrm{LC}}, \nabla^{\rm{Ch}}$ and $\nabla^{\rm{SB}}$
are all in $\sB_h$. We shall consider the preimage of them under the
isomorphism $\rho:\sA_{g,J}\>\sB_{h}$. By definition,
$\nabla^{\mathrm{LC}}\in \sA_{g,J}$ if and only if
$\nabla^{\rm{LC}}J=0$, i.e. $(M,g,J)$ is a K\"ahler manifold. It is
a natural question to find the preimage when it is not K\"ahler. The
following Lemma is well-known.

\blemma\label{canonicalfamily} Let $(M,g,J)$ be a Hermitian
manifold. Then

\bd \item $\rho^{-1}(\nabla^{\mathrm{Ch}})\in \sA_{g,J}$ is given by
\beq  g(\nabla _X Y,Z)=g(\LC_XY,Z)- \frac{1}{2}
d\om(JX,Y,Z).\label{realCh}\eeq

\item $\rho^{-1}(\hat\nabla^{\mathrm{LC}})\in \sA_{g,J}$ is given by \beq g(\nabla _X Y,Z)=g(\LC_XY,Z)+\frac{1}{4}
d\om(JX,JY,JZ)- \frac{1}{4} d\om(JX,Y,Z).\label{realLC} \eeq

\item $\rho^{-1}(\nabla^{\mathrm{SB}})\in \sA_{g,J}$ is given by \beq g(\nabla _X Y,Z)=g(\LC_XY,Z)+\frac{1}{2}
d\om(JX,JY,JZ).\label{realSB} \eeq

\ed \noindent Moreover, there exists a unique linear family
$\{\nabla^t\}_{t\in\R}\subset \sA_{g,J}$ such that
$$\rho(\nabla^0)=\nabla^{\mathrm{Ch}},\ \rho(\nabla^\frac{1}{2})=\hat\nabla^{\mathrm{LC}}
\qtq{and}\rho(\nabla^1)=\nabla^{\mathrm{SB}},$$ and it is given by
\beq g(\nabla ^{t} _X Y,Z)=g(\LC_XY,Z)+\frac{t}{2} d\om(JX,JY,JZ)+
\frac{t-1}{2} d\om(JX,Y,Z).\label{family} \eeq \elemma

\vskip 2\baselineskip

\section{Curvatures of connections in $\sA_{g,J}$.}
\label{curvatureofagj}

Recall that the real Chern connection
$\nabla^{\mathrm{Ch,\R}}=\rho^{-1}(\nabla^{\mathrm{Ch}})\in
\sA_{g,J}$ is given in (\ref{realCh}) and its complexificiation is
determined  by \beq
\begin{cases} \nabla^{\mathrm{Ch,\R}}_{\frac{\p}{\p
z^i}}\pzj=\Gamma_{ij}^k\pzk,\\
\nabla^{\mathrm{Ch,\R}}_{\frac{\p}{\p \bar z^i}}\pzj=0.
\end{cases}
\eeq For any $\nabla\in \sA_{g,J}$, there exists some
$A\in\Gamma(M,\Om^1(\mathrm{End(T_\R M)}))$ such that
$$\nabla=\nabla^{\mathrm{Ch,\R}}+A$$
and its complexficiation $\nabla:\Gamma(M,T_\C
M)\>\Gamma(M,\Om^1(T_\C M))$ is given by \beq \begin{cases}
\nabla^{}_{\frac{\p}{\p
z^i}}\pzj=\left(\Gamma_{ij}^k+\theta_{ij}^k\right)\pzk,\\
\nabla_{\frac{\p}{\p \bar z^i}}\pzj=-h_{j\bar q}h^{k\bar
p}\bar{\theta_{ip}^q}\pzk,
\end{cases} \qquad   \begin{cases}
\nabla^{}_{{\bpzi}}\bpzj=\left(\bar{\Gamma_{ij}^k}+\bar{\theta_{ij}^k}\right)\bpzk,\\
\nabla_{{\pzi}}\bpzj=-h_{q\bar j}h^{p\bar k}{\theta_{ip}^q}\bpzk.
\end{cases}
\label{theta}\eeq Moreover, its restriction to $T^{1,0}M$,
$\hat\nabla=\rho(\nabla)\in \sB_h$ is \beq
\begin{cases}
\hat \nabla^{}_{\frac{\p}{\p
z^i}}\pzj=\left(\Gamma_{ij}^k+\theta_{ij}^k\right)\pzk,\\
\hat\nabla_{\frac{\p}{\p \bar z^i}}\pzj=-h_{j\bar q}h^{k\bar
p}\bar{\theta_{ip}^q}\pzk.\label{hattheta}
\end{cases}
\eeq To make this correspondence more clear for the readers, we
write $\nabla^\theta$ for $\nabla\in\sA_{g,J}$ defined by
(\ref{theta}), and $\hat\nabla^\theta$ for $\hat\nabla\in\sB_h$
defined in (\ref{hattheta}).
 By using similar notations as in Section \ref{background}, the
 curvatures of $\nabla^\theta$ and $\hat\nabla^\theta$ are denoted by $R^\theta, \mathfrak{R}^\theta$. More precisely
\beq
R^{\theta}(X,Y,Z,W)=h(\nabla^\theta_X\nabla^\theta_YZ-\nabla^\theta_Y\nabla^\theta_XZ-\nabla^\theta_{[X,Y]}Z,W),\eeq
for  $X,Y,Z,W\in T_\C M$ and \beq
\ssR^{\theta}(X,Y,Z,W)=h(\hat\nabla^\theta_X\hat
\nabla^\theta_YZ-\hat \nabla^\theta_Y\hat \nabla^\theta_XZ-\hat
\nabla^\theta_{[X,Y]}Z,W),\eeq for  $ X,Y\in T_\C M, Z\in T^{1,0}M,
W\in T^{0,1}M$. We also use  conventions $R^\theta_{i\bar j
k\bar\ell}$ and $\ssR^\theta_{i\bar j k\bar \ell}$ for their
components.

\bproposition\label{fullcurvature0} For any $\nabla^\theta\in
\sA_{g,J}$ with
$\theta\in\Gamma\left(M,\Om^{1,0}\left(\mathrm{End}\left(T^{1,0}M\right)\right)\right)$,
the curvature tensors $R^\theta$ and $\ssR^\theta$ are determined by
\be\begin{cases} R^\theta_{i\bar jk\bar\ell}&=\ssR^\theta_{i\bar j
k\bar\ell} =\Ta_{i\bj
k\bel}-\left(h_{k\bbp}\frac{\p\bar{\ta_{j\ell}^p}}{\p
z^i}+h_{p\bel}\frac{\p\ta_{ik}^p}{\p
\bz^{j}}\right)+\left(\ta_{ik}^p\bar{\ta_{j\ell}^q}h_{p\bbq}-h^{m\bn}\ta_{im}^{p}\bar{\ta_{jn}^q}h_{p\bel}h_{k\bbq}\right),\\
R_{ijk\bar\ell}^\theta&=\ssR^\theta_{ijk\bel}=\left(\frac{\p\ta_{jk}^m}{\p
z^i}-\frac{\p\ta_{ik}^m}{\p
z^j}+\Gamma_{jk}^s\theta_{is}^{\ell}-\Gamma_{js}^{\ell}\theta_{ik}^s+\Gamma_{is}^{\ell}\theta_{jk}^s
-\Gamma_{ik}^s\theta_{js}^{\ell}\right)h_{m\bel}\\&+\left(\ta_{ip}^{m}\ta_{jk}^p-\ta_{ik}^p\ta_{jp}^{m}\right)h_{m\bel}.
 \end{cases}\ee

\eproposition

\bproof Since $\nabla^\theta J=0$, we have $$R^\theta_{i\bar
jk\bar\ell}=\ssR^\theta_{i\bar j k\bar\ell},\ \ \ R^\theta_{i
jk\bar\ell}=\ssR^\theta_{i j k\bar\ell}.  $$  Moreover, \be
\hat\nabla^\theta_\pzi\hat\nabla^\theta_{\pzj}\pzk
&=&\hat\nabla^\theta_\pzi\left(\Gamma_{jk}^\ell\pzl+\theta_{jk}^\ell\pzl\right)\\
&=&\left(\frac{\p\Gamma_{jk}^\ell}{\p
z^i}+\Gamma_{jk}^s(\Gamma_{is}^\ell+\theta_{is}^\ell)+\frac{\p\theta_{jk}^\ell}{\p
z^i}+\theta_{jk}^s(\Gamma_{is}^\ell+\theta_{is}^\ell)\right)\pzl,
\ee

\noindent and  \be (\ssR^\theta)_{ijk}^\ell=\left(\frac{\p
\theta_{jk}^{\ell}}{\p z^{i}}-\frac{\p \theta_{ik}^{\ell}}{\p z^{j}}
+\Gamma_{jk}^s\theta_{is}^{\ell}-\Gamma_{js}^{\ell}\theta_{ik}^s
+\Gamma_{is}^{\ell}\theta_{jk}^s-\Gamma_{ik}^s\theta_{js}^{\ell}\right)
+\left(\theta_{jk}^s\theta_{is}^\ell-\theta_{ik}^s\theta_{js}^\ell\right).\ee

\noindent Similarly, \be &&
h\left(\hat\nabla^\theta_\pzi\hat\nabla^\theta_{\bpzj}\pzk,\bpzl\right)\\
&=&\frac{\p}{\p
z^i}h\left(\hat\nabla^\theta_{\bpzj}\pzk,\bpzl\right)
-h\left(\hat\nabla^\theta_{\bpzj}\pzk,\hat\nabla^\theta_{\pzi}\bpzl\right)\\
&=& \frac{\p}{\p z^i}h\left
(-h_{k\bar{q}}h^{s\bar{p}}\bar{\theta_{jp}^q}\pzs,\bpzl\right)
-h\left(
h_{k\bar{p}}h^{s\bar{q}}\bar{\theta_{jq}^p}\pzs,h_{p\bar{\ell}}h^{q\bar{t}}\theta_{iq}^p\bpzt\right)
\\&=& -\frac{\p}{\p
z^i}\left(\bar{\theta_{j\ell}^s}h_{k\bar{s}}\right)
-h^{q\bar{p}}\theta_{iq}^{s}\bar{\theta_{jp}^t}h_{s\bar\ell}
h_{k\bar t}
\\&=& -\frac{\p}{\p
z^i}\bar{\theta_{j\ell}^s}h_{k\bar{s}}
 -\bar{\theta_{j\ell}^s}\Ga_{ik}^th_{t\bar s}
-h^{q\bar{p}}\theta_{iq}^{s}\bar{\theta_{jp}^t}h_{s\bar\ell}
h_{k\bar t},\ee

\noindent and \be &&
h\left(\hat\nabla^\theta_\bpzj\hat\nabla^\theta_{\pzi}\pzk,\bpzl\right)\\
&=&h\left(\hat\nabla^\theta_\bpzj\left((\Ga_{ik}^s+\ta_{ik}^s)\pzs\right),\bpzl\right)
\\&=&h\left(\frac{\p}{\p \bar{z}^j}
(\Ga_{ik}^s+\ta_{ik}^s)\pzs,\bpzl\right)+h\left((\Ga_{ik}^s+\ta_{ik}^s)
\left( -h_{s\bar p}h^{t\bar
q}\bar{\ta_{jq}^p}\right)\pzt,\bpzl\right)
\\&=&h_{s\bar\ell}\frac{\p \Ga_{ik}^s}{\p \bar{z}^j}
+h_{s\bar\ell}\frac{\p \ta_{ik}^s}{\p
\bar{z}^j}+(\Ga_{ik}^s+\ta_{ik}^s) ( -h_{s\bar
p}\bar{\ta_{j\ell}^p}). \ee Hence, we obtain the curvature formulas
in Proposition \ref{fullcurvature0}. \eproof

\noindent By using Proposition \ref{fullcurvature0}, one has
\bcorollary The first Ricci curvature of $\hat\nabla^\ta$ is \be
\ssRic^{(1)}(\ta)=\Ta^{(1)}-\I(\p\bar{\ta}_1-\bp\ta_1) \ee where
$\ta_1=\ta_{ik}^kdz^i$. \ecorollary

\noindent There are two important linear families  in $\sA_{g,J}$.
One
 is the Gauduchon family defined in (\ref{family}) and in this case,
\beq  \theta_{ij}^k=t\cdot T_{ij}^k\eeq and their curvatures are
given in Corollary \ref{god}. The other family is $\theta=t\cdot
\eta \ts \mathrm{Id_{T^{1,0}M}}$ for some one form $\eta=\eta_idz^i
\in\Gamma(M,\Om^{1,0}_M)$, and \beq \theta_{ij}^k=t\cdot\eta_i
\delta_{j}^k. \eeq

\bcorollary  The curvature formulas are \beq \ssR_{i\bj
k\bel}(\ta)=\Ta_{i\bj k\bel}-t\left(\frac{\p\bar{\eta_{j}}}{\p
z^i}+\frac{\p\eta_{i}}{\p \bz^{j}}\right)h_{k\bel},\ \ \
\ssR_{ijk\bel}(\ta)=t\left(\frac{\p\eta_{j}}{\p
z^i}-\frac{\p\eta_{i}}{\p z^j}\right)h_{k\bel},\eeq \beq
\ssRic^{(1)}(\ta)=\Ta^{(1)}-nt\I(\p\bar{\eta}-\bp\eta).\eeq
\ecorollary

\bremark When $d\eta=0$, one has
$$\ssR_{i\bj
k\bel}(\ta)=\Ta_{i\bj k\bel}, \qtq{and} \ssR_{ijk\bel}(\ta)=0$$ for
any $t\in \R$. \eremark

\vskip 2\baselineskip

\section{Geometry of real Chern-Einstein metrics}\label{cherneinstein}

In this section, we investigate real Chern-Einstein metrics and
prove Theorem \ref{main1}, Theorem \ref{main0} and  Theorem
\ref{main4}. Recall that the \emph{real Chern-Ricci curvature}
$\mathrm{Ric}(\nabla^{\mathrm{Ch,\R}},g)$ is defined by using the
Riemannian metric $g$, i.e.
$$\mathrm{Ric}(\nabla^{\mathrm{Ch,\R}},g)(X,Y)=\sum_{i=1}^{2n} R^{\mathrm{Ch,\R}}(X,e_i,e_i,Y)$$
where $\{e_i\}_{i=1}^{2n}$ is an orthonormal frame with respect to
$g$.

\bproposition The complexification of
$\mathrm{Ric}(\nabla^{\mathrm{Ch,\R}},g)$ is given by \beq
\mathrm{Ric}(\nabla^{\mathrm{Ch,\R}},g)=\Theta^{(3)}_{i\bar
j}dz^i\ts d\bar z^j+\Theta^{(4)}_{i\bar j}d\bar z^j\ts dz^i
\label{complexricci}\eeq where $\Theta_{ij}^{(3)}$ and
$\Theta_{ij}^{(4)}$ are defined in (\ref{34}). \eproposition

\bproof By using Theorem \ref{main2} for $\theta=0$, we have \be
\mathrm{Ric}(\nabla^{\mathrm{Ch,\R}},g)\left(\pzi,\bpzj\right)&=&
h^{k\bar\ell }R_{ki \bar j \bar \ell}+h^{\bar\ell k}{R}_{\bar \ell
i\bar j k }\\
&=&h^{k\bar\ell} R_{i\bar\ell k\bar
j}=h^{k\bar\ell}\Theta_{i\bar\ell k\bar j}\\
&=&\Theta^{(3)}_{i\bar j}\ee where $R$ stands for
$R^{\mathrm{Ch,\R}}$. Similarly, we have
$$\mathrm{Ric}(\nabla^{\mathrm{Ch,\R}},g)\left(\bpzj,\pzi\right)=\Theta^{(4)}_{i\bar
j}$$ and
$\mathrm{Ric}(\nabla^{\mathrm{Ch,\R}},g)\left(\pzi,\pzj\right)=
\mathrm{Ric}(\nabla^{\mathrm{Ch,\R}},g)\left(\bpzj,\bpzi\right)=0$.
\eproof

\bdefinition $(M,g,J,\nabla^{\mathrm{Ch,\R}})$ is called \emph{real
Chern-Einstein} if \beq
\mathrm{Ric}(\nabla^{\mathrm{Ch,\R}},g)=\lambda g \qtq{for some
$\lambda\in \R.$}\eeq  If
$\mathrm{Ric}(\nabla^{\mathrm{Ch,\R}},g)=0$, it is also called
\emph{real Chern-Ricci flat}. Moreover,
$(\nabla^{\mathrm{Ch,\R}},g)$ has positive  real Chern-Ricci
curvature if
$\mathrm{Ric}(\nabla^{\mathrm{Ch,\R}},g)\in\Gamma(M,\mathrm{Sym}^{\ts
2}T_\R M)$ and it is positive  definite. The negativity can be
defined similarly. \edefinition

\btheorem\label{positive} Let $(M,g,J)$ be a Hermitian manifold.
Then $\mathrm{Ric}(\nabla^{\mathrm{Ch,\R}},g)$ is positive if and
only if \beq \Theta^{(1)}-\p\p^*\omega \eeq is a positive definite
Hermitian $(1,1)$ form. In particular, $(\nabla^{\mathrm{Ch,\R}},g)$
is real Chern-Einstein with constant $\lam\in \R$ if and only if
\beq \Ta^{(1)}-\p\p^\st\om=\lam\om. \label{realeinstein}\eeq where
$\Theta^{(1)}$ is the first Chern-Ricci curvature.
 \etheorem

\bproof By formula (\ref{complexricci}), if
$\mathrm{Ric}(\nabla^{\mathrm{Ch,\R}},g)\in\Gamma(M,\mathrm{Sym}^{\ts
2}T_\R M)$, then $\Theta^{(3)}_{i\bar j}=\Theta^{(4)}_{i\bar j}$.
Hence,
$$\mathrm{Ric}(\nabla^{\mathrm{Ch,\R}},g)=\Theta^{(3)}_{i\bar
j}(dz^i\ts d\bar z^j+d\bar z^j\ts dz^i).$$ Therefore, the real
Chern-Ricci curvature is positive definite if and only if
$\Theta^{(3)}_{i\bar j}v^i\bar v^j>0$ for every nonzero vector
$(v^i)$. By Theorem \ref{chernricci} or Corollary \ref{god} (when
$t=0$), we know $$\Theta^{(3)}=\sq \Theta^{(3)}_{i\bar j}dz^i\wedge
d\bar z^j=\Theta^{(1)}-\p\p^*\omega$$ is a positive definite
Hermitian $(1,1)$ form. In particular, $(\nabla^{\mathrm{Ch,\R}},g)$
is real Chern-Einstein with constant $\lam\in \R$ if and only if
(\ref{realeinstein}) holds.
 \eproof

\noindent\emph{Proof of Theorem \ref{main1}.} By applying $\p$ to
the equation (\ref{realeinstein}), we have $\lambda \p\omega=0$.
Hence if $\lambda\neq 0$, $d\omega=0$ and $(M,g,J)$ is
K\"ahler-Einstein.\qed \\

\noindent \emph{Proof of Theorem \ref{main0}.} Suppose
$\mathrm{Ric}(\nabla^{\mathrm{Ch,\R}},g)>0$, by Theorem
\ref{positive} we deduce the Hermitian $(1,1)$-form
$$\Om_0:=\Theta^{(1)}-\p\p^*\omega>0$$
and $\p\Om_0=0$. Since $\Omega_0$ is also real, we obtain
$d\Om_0=0$. Hence, $\Om_0$ is a K\"ahler form and $(M,J)$ is a
K\"ahler manifold. Moreover, the $(1,1)$-form $\p\p^*\omega$ is both
$d$-closed and $\p$-exact. By $\p\bp$-lemma on the K\"ahler manifold
$(M,J)$, there exists some $f\in C^\infty(M,\R)$ such that
$\p\p^*\omega=-\sq\p\bp f$. Hence $\Theta^{(1)}=\Om_0+\sq \p\bp f$
and so $c_1(M,J)>0$. The proof for
$\mathrm{Ric}(\nabla^{\mathrm{Ch,\R}},g)<0$ is similar.
\qed\\

When $\lambda=0$, we have the following result.
\bcorollary\label{realchernricciflat} Let $(M,g,J)$ be a Hermitian
manifold.  Then $(\nabla^{\mathrm{Ch,\R}},g)$ is real Chern-Ricci
flat if and only if
$$\Theta^{(1)}-\p\p^*\omega=\Theta^{(1)}-\bp\bp^*\omega=0.$$
\ecorollary

\noindent\emph{Proof of Theorem \ref{main4}}. On the standard Hopf
manifold $\S^{2n-1}\times \S^1$ with canonical metric $$\omega_0=\sq
h_{i\bar j}dz^i\wedge d\bar z^j= \frac{4\delta_{i\bar j}}{|z|^2}\sq
dz^i\wedge  d\bar z^j$$ we know the following metric
  \beq \omega=\omega_0-\frac{4}{n}\cdot
\sq\p\bp\log|z|^2\eeq is Levi-Civita-Ricci flat
(\cite[Theorem~6.2]{LY16} or Theorem \ref{main6} with
$t=\frac{1}{2}$). On the other hand, one can show directly (see also
\cite[Theorem~6.2]{LY16}) that
$$\p\p^* \omega=n\sq\p\bp\log |z|^2$$
where $\p^*$ is taken with respect to $\omega$. By Corollary
\ref{realchernricciflat}, we deduce that $(\nabla^{\mathrm{Ch,\R}},
\omega)$ is real Chern-Ricci flat.\qed\\

\btheorem\label{main7} Let $(M,g,J)$ be a compact Hermitian
manifold. Suppose $(\nabla^{\mathrm{Ch,\R}},g)$ is real Chern-Ricci
flat. Then  $(M,J)$ is one of the following
 \begin{enumerate}
 \item  $(M,J)$ is K\"ahler:  $c_1(M,J)=0$, i.e. $(M,J)$ has a
 K\"ahler-Ricci flat metric;

\item $(M,J)$ is not K\"ahler:  \bd \item $\kappa(M)=0$ and $c_1^{\mathrm{BC}}(M)=0$. Moreover,  $(M,J)$ has a balanced metric and $K_M$  holomorphic torsion, i.e. $K_M^{\otimes m}=0$ for some $m\in \N_+$.

\item  $\kappa(M)=-\infty $ and $c_1^{\mathrm{AC}}(M)=0$.

\ed
 \end{enumerate}

\noindent Moreover, the only non-K\"ahler compact complex surface
which can support real Chern-Ricci metric is the Hopf surface.

\etheorem

\noindent \emph{Proof of Theorem \ref{main7}.}  By Theorem
\ref{positive},  we have $\Theta^{(1)}=\p\p^*\omega$. Hence, the
real $(1,1)$ form $\p\p^*\omega$ is $d$-closed and $\p$-exact.

 If
$(M,J)$ is K\"ahler, by $\p\bp$-lemma, there exists some $f\in
C^\infty(M,\R)$ such that $\Theta^{(1)}=\p\p^*\omega=\sq\p\bp f$. In
particular, $c_1(M)=0$ on the K\"ahler manifold $(M,J)$. By the
Calabi-Yau theorem, there exists a  K\"ahler-Ricci flat metric
$\tilde \omega$ which is possibly different from $\omega$.

When $(M,J)$ is not K\"ahler, by Theorem \cite[Theorem~5.2]{LY17} or
Theorem \ref{kod} ($t=\frac{1}{2}$), we obtain the classification.
When $\dim_\C M=2$, by Theorem \ref{cgrfs} ($t=\frac{1}{2}$), we
know the Hopf surface is the only non-K\"ahler surface which can
support real Chern-Ricci flat metrics. \qed

\bremark\label{remark} On a K\"ahler Calabi-Yau manifold $M$, there
exist non-K\"ahler metrics which are real Chern-Ricci flat. Indeed,
let $\omega_{\mathrm{CY}}$ be a Calabi-Yau K\"ahler metric on $M$.
 Then for any non constant smooth function $f\in C^\infty(M,\R)$,   by Yau's  theorem (\cite{Yau78}), there
 exists a K\"ahler metric $\omega_0$ such that
 $$\omega_0^n=e^{-f}\omega^n_{\mathrm{CY}}.$$
Let $\omega_f=e^f\omega_0$. We have
$$\bp^*_f\omega_f=\bp^*_0\omega_0+(n-1)\sq\p f=(n-1)\sq\p f.$$
Hence $\bp\bp_f^*\omega_f=\p\p^*_f\omega_f=-\sq(n-1)\p\bp f$.
Moreover, we have
$$\omega_f^{n}=e^{(n-1)f}\omega^n_{\mathrm{CY}}$$ which implies
$\Theta^{(1)}(\omega_f)=\Theta^{(1)}(\omega_{\mathrm{CY}})-(n-1)\sq\p\bp
f=\p\p^*_f\omega_f$. By Corollary \ref{realchernricciflat},
$\omega_f$ is a real Chern Ricci-flat metric, and it is a
non-K\"ahler metric.
 \eremark

\vskip 2\baselineskip

\section{Classification of compact complex surfaces with $t$-Gauduchon-Ricci flat
metrics}\label{classification}

 In this section, we
classify compact complex surfaces with $t$-Gauduchon-Ricci flat
metrics.  One of the key ingredients is understanding the geometry
of scalar curvatures of Gauduchon connections. The following theorem
generalizes results in \cite{LY16, LY17, Yang19,
 HLY18} to the Gauduchon family.

  \btheorem\label{kod}
 Let $M$ be a compact complex manifold. Suppose $\om$ is a Hermtian metric and $\nabla^t$ is the Gauduchon connection of $M$. If $s^{(1)}(\omega, t)\geq 0$ for some $t>0$, then either
 \begin{enumerate}
 \item $\kappa(M)=-\infty $; or
 \item $ \kappa(M)=0$ and $(M,\om)$ is conformally balanced with $K_M$  holomorphic torsion, i.e. $K_M^{\otimes m}=0$ for some $m\in \N_+$.
 \end{enumerate}
 \etheorem

 \noindent  For  $t<0$, we have a similar result that
  \btheorem\label{tneg}
 Let $M$ be a compact complex manifold.
 Suppose $\om$ is a Hermtian metric and $\nabla^t$ is the Gauduchon connection of $M$. If $s^{(1)}(\omega, t)\leq 0$ for some $t<0$
 and $s^{(1)}(\omega, t)$ is strictly negative at some point, then $K^{-1}_M $ is not pseudo-effective.
 \etheorem

To prove the Theorem \ref{kod} and Theorem \ref{tneg}, we first
calculate the total scalar curvature of the  Gauduchon metric in the
conformal class of $\om$. It is well-known that there exists a
smooth  function $f$ on $M$ such that $\om_f=e^f\om$ is Gauduchon,
i.e. $\p\bp\om_{f}^{n-1}=0$.

\blemma\label{gts} Let $s_f$ be the Chern scalar curvature of
$\om_f=e^f\om$. Then we  have: \be\label{gtseqn} \int_M s_f
\frac{\om_f^n}{n!}= \int_M f^{n-1}\cdot s^{(1)}(\omega,
t)\cdot\frac{\om^n}{n!} +t\int_M \left(|\bp_f^\st
\om_f|_f^2+|\p_f^\st \om_f|_f^2\right)\frac{\om_f^{n}}{n!}.\ee
\elemma \bpf We have the relation
$\Theta^{(1)}(\om_f)=\Theta^{(1)}(\om)-\I n\p\bp f$, it then follows

\be \int_ M s_f  \frac{\om_f^n}{n!}= \int_M \Theta^{(1)}(\om_f)\w
\frac{\om_f^{n-1}}{(n-1)!}&=&\int_M (\Theta^{(1)}(\om)-n\sq\p\bp
f)\w \frac{\om_f^{n-1}}{(n-1)!}
\\&=& \int_M \Theta^{(1)}(\om)\w \frac{\om_f^{n-1}}{(n-1)!}.
\ee

\noindent By using Corollary \ref{god},
$\ssRic^{(1)}(\om,t)=\Theta^{(1)}(\om)-t(\bp\bp^\st\om+\p\p^\st\om)
$,
 \be
\int_M s_f  \frac{\om_f^n}{n!} &=& \int_M \Theta^{(1)}(\om)\w
\frac{\om_f^{n-1}}{(n-1)!}
\\&=&\int_M (\ssRic^{(1)}(\om,t)+t\bp\bp^\st\om+t\p\p^\st\om)\w \frac{\om_f^{n-1}}{(n-1)!}
\\&=& \int_M e^{(n-1)f}\cdot s^{(1)}(\omega, t)\cdot \frac{\om^n}{n!}+t\int_M (\bp\bp^\st\om+\p\p^\st\om)\w \frac{\om_f^{n-1}}{(n-1)!}.
\ee

\noindent By using the formula  \cite[Lemma~3.4]{Yang20} \be
\bp_f^\st\om_f=\bp^\st\om+(n-1)\I \p f, \ee we get \be \bp
\bp_f^\st\om_f=\bp\bp^\st\om-(n-1)\I \p \bp f,\ \ \ \ \p
\p_f^\st\om_f=\p\p^\st\om-(n-1)\I \p \bp f.\ee

\noindent Therefore, \be \int_M (\bp\bp^\st\om+\p\p^\st\om)\w
\frac{\om_f^{n-1}}{(n-1)!} &=&\int_M
\left(\bp_f\bp_f^\st\om+\p_f\p_f^\st\om\right)\w
\frac{\om_f^{n-1}}{(n-1)!}
\\&=& \int_M \left(|\bp_f^\st \om_f|_f^2+|\p_f^\st \om_f|_f^2\right)\frac{\om_f^{n}}{n!} .
\ee Hence, (\ref{gtseqn}) follows. \epf

As an application of Lemma \ref{gts}, we can prove Theorem \ref{kod}
and \ref{tneg}.

\bpf[Proof of Theorem \ref{kod}] If $s^{(1)}(\omega, t)\geq0$ for
some $t>0$, by Lemma \ref{gts}, we have \be \int_M s_f
\frac{\om_f^n}{n!}\geq   t\int_M \left(|\bp_f^\st
\om_f|_f^2+|\p_f^\st \om_f|_f^2\right)\frac{\om_f^{n}}{n!}\geq 0.\ee

\noindent If  $\int_ M s_f  \frac{\om_f^n}{n!} >0$, by
\cite[Corollary~3.3]{Yang19}, we have $\kappa(M)=-\infty$. If $\int_
M s_f \frac{\om_f^n}{n!}=0$, then we must have $\bp^\st\om_f=0$,
i.e. $\om_f$ is a balanced metric. Hence $\om$ is conformal
balanced. In this case, also by \cite[Theorem~1.4]{Yang19}, it
follows either $\kappa(M)=-\infty$ or $\kappa(M)=0$ with $K_M$  a
holomorphic torsion. \epf

\bpf[Proof of Theorem \ref{tneg}] By using Lemma \ref{gts} again, we
deduce
 \be
\int_M s_f  \frac{\om_f^n}{n!}\leq   \int_M e^{(n-1)f}\cdot
s^{(1)}(\omega, t) \frac{\om^n}{n!}< 0. \ee
 By \cite[Theorem~1.3]{Yang19}, $K_M^{-1}$ is not pseudo-effective.
 \epf

\vskip 1\baselineskip

\noindent Now we are ready to establish the classification.
\btheorem\label{cgrfs} Let $S$ be a compact complex surface. If it
admits a $t$-Gauchon-Ricci flat metric $\omega$ for some $t>0$, then
$S$ is a minimal surface lying in one of the following:
\begin{enumerate}
\item Enriques surfaces.
\item bi-elliptic surfaces.
\item K3-surfaces.
\item  2-tori.
\item Hopf surfaces.
\end{enumerate}
\etheorem

\noindent We shall prove Theorem \ref{cgrfs} following ideas in
\cite{HLY18}. By Theorem \ref{kod}, we have: \bcorollary Suppose
$t>0$. If a complex surface $S$ can admit a $t$-Gaudchon-Ricci flat
metric, then
\begin{enumerate}\label{grfp1}
\item either $\kappa(S)=-\infty$ ; or
\item $\kappa(S)=0$. In this case, $\om$ is conformal K\"ahler and $K_S$ is a holomorphic torsion, i.e. $K_S^{\otimes m}=0$ for some integer $m\in \Z$.
\end{enumerate}
In both cases, we have $c_1^2(S)=0$. \ecorollary

\noindent  We need two more lemmas and their proofs are similar to
those in \cite{HLY18}. \blemma\cite[Theorem~4.3]{HLY18}\label{grfp4}
Let $S$ be a complex surface with $\kappa(S)=-\infty$. If $S$ admits
a $t$-Gauduchon-Ricci flat metric, then $S$ must be non-K\"ahler.
\elemma

 \blemma\cite[Theorem~5.1]{HLY18}\label{grfp2} Let $S$ be a
non-K\"ahler complex surface with $\kappa(S)=-\infty$. If
$c_1^2(S)=0$, then $S$ must be minimal. \elemma

\noindent As an application of Proposition \ref{grfp1}, Lemma
\ref{grfp4} and Lemma \ref{grfp2}, one has \bcorollary\label{grfp3}
If a compact complex surface $S$ admits a $t$-Gauduchon-Ricci flat
metric $\om$, then $S$ must be minimal. \ecorollary

\vskip 1\baselineskip

\noindent \emph{Proof of Theorem \ref{cgrfs}.}  Supposee $S$
supports a $t$-Gauduchon-Ricci flat metric, then by Proposition
\ref{grfp1}, $\kappa(S)\leq 0$ and by Proposition \ref{grfp3}, $S$
is also minimal.\\

\noindent (A). $\kappa(S)=0$. By Kodaira-Enriques' classification,
$S$ is exactly one of the following:
\begin{enumerate}
\item Enrique surface;
\item bielliptic surface;
\item K3 surface;
\item torus.
\end{enumerate}
They are all K\"ahler Calabi-Yau.\\

\noindent (B).  $\kappa(S)=-\infty$. By using Kodaira-Enriques'
classification again, $S$ can only be one of the following:
\begin{enumerate}
\item minimal rational surface;
\item ruled surface of $g>0$;
\item surface of class $\text{VII}_0$.
\end{enumerate}
By proposition \ref{grfp4}, $S$ is non-K\"ahler and so $S$ can only
be a surface of class $\text{VII}_0$:
\begin{enumerate}
\item class $\text{VII}_0$ surfaces with type $b_2>0$.
\item Inoue surface.
\item Hopf surface.
\end{enumerate}
By using similar strategies as in the proof of \cite[Theorem
5.1]{HLY18}, one can show $S$ can only be a Hopf surface. \qed

 \bremark An explicit
$t$-Gauduchon-Ricci flat metrics on diagonal Hopf surface is
constructed in Theorem \ref{ricciflat}. \eremark

\vskip 2\baselineskip

\section{Explicit construction of $t$-Gauduchon-Ricci flat  metrics on Hopf
manifolds} \label{hopfmanifold}
 Let $M=\S^{2n-1}\times \S^1$ be the standard
$n$-dimensional ($n\geq 2$) Hopf manifold. It is diffeomorphic to
$\C^n- \{0\}/G$ where $G$ is a cyclic group generated by the
transformation $z\rightarrow \frac{1}{2}z$. It has an induced
complex structure from $\C^n-\{0\}$.   On $M$, there is a natural
induced metric $\omega_0$ given by
 \beq \omega_0=\sq h_{i\bar j}dz^i\wedge d\bar z^j=\sq
\frac{4\delta_{i\bar j}}{|z|^2}dz^i\wedge  d\bar z^j.\label{hopf}
\eeq \noindent The main result of this section is the following:

\btheorem\label{ricciflat} The Hermitian metric \beq
\Om_t=\omega_0+4\left(\frac{2(n-1)t}{n}-1\right)\cdot\sq \p\bp\log
|z|^2,\ \ t>0 \eeq is $t$-Gauduchon-Ricci flat, i.e.,
$\ssRic^{(1)}(\Omega_t, t)=0$. \etheorem

\bproof We shall use similar constructions as in
\cite[Section~6]{LY16}. More precisely, we consider the perturbed
Hermitian metric \beq \omega_\lambda=\omega_0+4\lambda\sq\p\bp\log
|z|^2,\qtq{with} \lambda>-1.\eeq

\noindent It is shown in \cite[Theorem~6.2]{LY16} that \beq
\Theta^{(1)}(\omega_\lambda)=n\cdot \sq\p\bp\log |z|^2\eeq and \beq
\frac{\p\p^*\omega_\lambda+\bp\bp^*\omega_\lambda}{2}=
\frac{n-1}{1+\lambda}\cdot\sq\p\bp\log |z|^2. \eeq By using
Corollary \ref{god}, we deduce \beq \ssRic^{(1)}(\omega_\lambda,
t)=\Theta^{(1)}-t(\p\p^*\omega_\lambda+\bp\bp^*\omega_\lambda)=\left(n-\frac{2(n-1)t}{1+\lambda}\right)\sq
\p\bp\log |z|^2.\eeq Therefore,
$\mathrm{Ric}^{(1)}(\omega_\lambda,t)=0$ if and only if
$\displaystyle{\lambda=\frac{2(n-1)t}{n}-1}$. \eproof


\bremark Note that when $t=0$, $\Om_0=\omega_0-4\sq\p\bp\log|z|^2$
is not a Hermitian metric since the corresponding matrix is not
positive definite. It is also well-known that there is no
Chern-Ricci flat metrics on $\S^{2n-1}\times \S^1$, although there
are Levi-Civita-Ricci flat metrics $(t=\frac{1}{2})$ and
Strominger-Bismut-Ricci flat metrics $(t=1)$. On the Hopf surface
$\S^3\times \S^1$, the canonical metric $\Om_1=\omega_0$ is
Strominger-Bismut-Ricci flat.

 \eremark

\vskip 2\baselineskip

\section*{Appendix: curvatures of Gauduchon
connections}\label{appendix}

In this section, we give a detailed proof of Lemma
\ref{canonicalfamily} (which is definitely well-known to experts)
and  Corollary \ref{god}.

Let $\{\nabla^{\lambda,\mu}\}_{\lambda,\mu\in \R}\subset \sA_g$ be a
family affine connections on  $T_\R M$ defined by \beq g(\nabla
^{\lambda,\mu} _X Y,Z):=g(\LC_XY,Z)+\lambda d\om(JX,JY,JZ)+\mu
d\om(JX,Y,Z), \eeq for $X,Y,Z\in \Gamma(M,T_\R M)$. Let $\{z^i\}$ be
the local holomorphic coordinates on $M$. We consider the
complexification of $\nabla^{\lambda,\mu}$ by setting \beq
\displaystyle{\begin{cases}
\nabla^{\lambda,\mu}_{\pzi}\pzj=\Gamma_{ij}^k(\lambda,\mu)\pzk+\Gamma_{ij}^{\bk}(\lambda,\mu)\bpzk,\\
\cd^{\lam,\mu}_{\bpzi}\pzj=\Ga_{\bi j}^k(\lam,\mu)\pzk+\Ga_{\bi
j}^{\bk}(\lam,\mu)\bpzk. \end{cases} }\eeq

\blemma\label{fullchristoffel} We have the following relations
\be\label{Crstf-eqn1}
\Ga_{ij}^k(\lam,\mu)=\Gamma_{ij}^k-\left(\lam+\mu+\half
\right)T_{ij}^k, \quad \Ga_{ij}^{\bk}(\lam,\mu)=0; \ee and
\be\label{Crstf-eqn2} \Ga^{k}_{\bi
j}(\lam,\mu)=\left(\lam+\mu+\half\right)h^{k\bm}h_{j\bn}\bar{T_{im}^n},
\quad \Ga_{\bi
j}^{\bk}(\lam,\mu)=\left(-\lam+\mu+\half\right)h^{m\bk}h_{n\bi}T^n_{jm}
\ee where $\Gamma_{ij}^k=h^{k\bar \ell}\frac{\p h_{j\bar\ell}}{\p
z^i}$ and $T_{ij}^k=\Gamma_{ij}^k-\Gamma_{ji}^k$. \elemma

\bproof It follows from standard computations, and for readers'
convenience  we include a straightforward proof here. At first, we
have \beq \LC_{\frac{\p}{\p z^i}}\frac{\p}{\p
z^j}=\frac{1}{2}h^{k\bar \ell}\left(\frac{\p h_{j\bar \ell}}{\p
z^i}+\frac{\p h_{i\bar \ell}}{\p z^j}\right) \frac{\p}{\p z^k} \eeq
and \beq \LC_{\frac{\p}{\p \bar z^i}}\frac{\p}{\p z^j}=\frac{1}{2}
h^{k\bar \ell}\left( \frac{\p h_{j\bar\ell}}{\p\bar z^i}-\frac{\p
h_{j\bar i}}{\p\bar z^\ell}\right)\frac{\p}{\p
z^k}+\frac{1}{2}h^{k\bar q}\left(\frac{\p h_{k\bar i}}{\p
z^j}-\frac{\p h_{j\bar i}}{\p z^k}\right)\frac{\p}{\p \bar z^q}.
\eeq Note also that
$d\om\left(\pzi,\pzj,\bpzl\right)=\sq\left(\frac{\p h_{j\bar
\ell}}{\p z^i}-\frac{\p h_{i\bar \ell}}{\p z^j}\right)$.  Hence, we
have \be &&h\left(\nabla^{\lambda,\mu}_{\pzi}\pzj,\bpzl\right)
\\&=&h\left(\LC_{\pzi}\pzj,\bpzl\right)+\lam d\om\left(J{\pzi},J\pzj,J\bpzl\right)+\mu d\om\left(J{\pzi},\pzj,\bpzl\right)
\\&=&\frac{1}{2}\left(\frac{\p h_{j\bar \ell}}{\p
z^i}+\frac{\p h_{i\bar \ell}}{\p z^j}\right) +\lam \I
d\om\left(\pzi,\pzj,\bpzl\right)+\I \mu
d\om\left(\pzi,\pzj,\bpzl\right)
\\&=&\frac{\p h_{j\bar\ell}}{\p z^i}-\left(\lam+\mu+\half
\right)\left(\frac{\p h_{j\bar\ell}}{\p z^i}-\frac{\p
h_{i\bar\ell}}{\p z^j}\right). \ee Therefore,
$$\Ga_{ij}^k(\lam,\mu)=\Gamma_{ij}^k-\left(\lam+\mu+\half
\right)T_{ij}^k.$$

\noindent Similarly, one has \be
&&h\left(\nabla^{\lambda,\mu}_{\pzi}\pzj,\pzk\right)
\\&=&h\left(\LC_{\pzi}\pzj,\pzk\right)+\lam d\om\left(J{\pzi},J\pzj,J\pzk\right)+\mu d\om\left(J{\pzi},\pzj,\pzk\right)
\\&=&0 \ee
since the metric is $J$-invariant. Therefore,
$\Ga_{ij}^{\bk}(\lam,\mu)=0$.

For the second part, we have $d\om\left(\bpzi,\pzj,\bpzl\right)=\sq
\left(\frac{\p h_{j\bar \ell}}{\p\bar z^i}-\frac{\p h_{j\bar
i}}{\p\bar z^\ell}\right)$ and
 \be
&&h\left(\nabla^{\lambda,\mu}_{\bpzi}\pzj,\bpzl\right)
\\&=&h\left(\LC_{\bpzi}\pzj,\bpzl\right)+\lam d\om\left(J{\bpzi},J\pzj,J\bpzl\right)+\mu d\om\left(J{\bpzi},\pzj,\bpzl\right)
\\&=&\frac{1}{2}\left(\frac{\p h_{j\bar \ell}}{\p
\bar z^i}-\frac{\p h_{j\bar i}}{\p \bar z^\ell}\right) -\lam \I
d\om\left(\bpzi,\pzj,\bpzl\right)-\I \mu
d\om\left(\bpzi,\pzj,\bpzl\right)
\\&=&\frac{1}{2}\left(\frac{\p h_{j\bar \ell}}{\p
\bar z^i}-\frac{\p h_{j\bar i}}{\p \bar
z^\ell}\right)+(\lambda+\mu)\left(\frac{\p h_{\ell\bar j}}{\p\bar
z^i}-\frac{\p h_{j\bar i}}{\p\bar
z^\ell}\right)=\left(\lambda+\mu+\frac{1}{2}\right)\left(\frac{\p
h_{\ell\bar j}}{\p\bar z^i}-\frac{\p h_{j\bar i}}{\p\bar
z^\ell}\right). \ee Therefore,
$$\Ga^{k}_{\bi
j}(\lam,\mu)=\left(\lam+\mu+\half\right)h^{k\bm}h_{j\bn}\bar{T_{im}^n}.$$

\noindent Similarly, one has \be
&&h\left(\nabla^{\lambda,\mu}_{\bpzi}\pzj,\pzk\right)
\\&=&h\left(\LC_{\bpzi}\pzj,\pzk\right)+\lam d\om\left(J{\bpzi},J\pzj,J\pzk\right)+\mu d\om\left(J{\bpzi},\pzj,\pzk\right)
\\&=&\frac{1}{2}\left(\frac{\p h_{k\bar i}}{\p
z^j}-\frac{\p h_{j\bar i}}{\p
z^k}\right)+\sq(\lambda-\mu)d\om\left(\bpzi,\pzj,\pzk\right)\\
&=& \left(-\lambda+\mu+\frac{1}{2}\right)\left(\frac{\p h_{k\bar
i}}{\p z^j}-\frac{\p h_{j\bar i}}{\p z^k}\right)\ee and we deduce
$$\Ga_{\bi
j}^{\bk}(\lam,\mu)=\left(-\lam+\mu+\half\right)h^{m\bk}h_{n\bi}T^n_{jm}.$$
The proof of Lemma \ref{fullchristoffel} is completed.
 \eproof

\noindent\emph{Proof of Lemma \ref{canonicalfamily}.} By Lemma
\ref{fullchristoffel} and formula (\ref{realcompatible}), we deduced
that $\nabla^{\lambda,\mu}\in \sA_{g,J}$ if and only if $\Ga_{\bi
j}^{\bk}(\lam,\mu)=0$, i.e. \beq \nabla^{\lambda,\mu}\in
\sA_{g,J}\Longleftrightarrow
\left(-\lambda+\mu+\frac{1}{2}\right)d\omega=0.\label{Chern0}\eeq

\noindent On the other hand,  the restricted connection $\hat
\nabla^{\lambda,\mu}=\pi(\nabla^{\lambda,\mu})$ on the holomorphic
tangent bundle $T^{1,0}M$ is determined by \beq \hat
\nabla^{\lambda,\mu}_{\pzi}\pzj=\Gamma_{ij}^k(\lambda,\mu)\pzk,\ \
\hat \nabla^{\lambda,\mu}_{\bpzi}\pzj=\Ga_{\bi j}^k(\lam,\mu)\pzk
\eeq where \be
\Ga_{ij}^k(\lam,\mu)=\Gamma_{ij}^k-\left(\lam+\mu+\half
\right)T_{ij}^k,\ \ \Ga^{k}_{\bi
j}(\lam,\mu)=\left(\lam+\mu+\half\right)h^{k\bm}h_{j\bn}\bar{T_{im}^n}.\ee

 \noindent  Recall that the Chern connection $\nabla^{\mathrm{Ch}}$
of $T^{1,0}M$ is characterized by
$$\nabla^{\mathrm{Ch}}_{\pzi}\pzj=\Gamma_{ij}^k\pzk,\ \
\ \ \nabla^{\mathrm{Ch}}_{\bpzi}\pzj=0. $$ Hence,  \beq
\pi(\nabla^{\lambda,\mu})=\nabla^{\mathrm{Ch}}\Longleftrightarrow
\left(\lambda+\mu+\frac{1}{2}\right)d\omega=0.\label{Chern1}\eeq By
using (\ref{Chern0}) and (\ref{Chern1}), we deduce \beq
\rho(\nabla^{\lambda,\mu})=\nabla^{\mathrm{Ch}}\Longleftrightarrow
d\omega=0, \qtq{or} (\lambda,\mu)=\left(0,-\frac{1}{2}\right).\eeq
Thus, we obtain (\ref{realCh}). Similarly, one can show
(\ref{realLC}) and (\ref{realSB}).  The uniqueness of the family
(\ref{family}) follows from the linear property.\qed

\subsection{Curvature formula of Gauduchon connections} Recall that
there is a linear family of connections  defined by
$$g(\nabla ^{t} _X Y,Z)=g(\LC_XY,Z)+\frac{t}{2}
d\om(JX,JY,JZ)+ \frac{t-1}{2} d\om(JX,Y,Z).$$ By using Lemma
\ref{fullchristoffel}, the Gauduchon connection is determined by:
\beq \dt_{\pzi}\pzj=\Ga_{ij}^k(t)\pzk \quad\text{and} \quad
\dt_{\bpzi}\pzj=\Ga_{\bi j}^k(t)\pzk, \eeq where the coefficients
$\Ga_{ij}^k$ and $\Ga_{\bi j}^k$ are given as follows:
\be\label{Crstf-eqn3} \Ga_{ij}^k(t)=\Ga_{ij}^k-t T_{ij}^k\quad
\text{and} \quad \Ga_{\bi j}^k(t)=t\cdot
h^{k\bm}h_{j\bn}\bar{T_{im}^n}. \ee

\btheorem\label{fullcurvature}
 The curvature tensor of Gauduchon connection $\nabla^t$ is given by:
\beq R_{i\bj k\bel}(t) =\Ta_{i\bj k\bel}+t(\Ta_{i\bel
k\bj}+\Ta_{k\bj i\bel}-2\Ta_{i\bj k\bel})+t^2(T_{ik}^p\bar{T_{j\ell
}^q}h_{p\bbq}-h^{p\bbq}h_{m\bar\ell} h_{k\bar
n}T_{ip}^m\bar{T_{jq}^n}).\label{curvature} \eeq
 \etheorem

\bproof In setting of Proposition \ref{fullcurvature0}, we have
$\theta_{ij}^k=-tT_{ij}^k$. Hence (\ref{curvature}) follows from
Proposition \ref{fullcurvature0} and the relation $\displaystyle{
\frac{\p T_{ik}^{\ell}}{\p\bar z^j}=-\Ta_{i\bj k}^{\ell}+\Ta_{k\bj
i}^{\ell}}.$ \eproof

\noindent \emph{Proof of Corollary \ref{god}.} It follows from
Theorem \ref{chernricci} and Theorem \ref{fullcurvature}.

\vskip 2\baselineskip

\end{document}